\documentclass[12pt]{amsart}
\usepackage{color}
\newcommand{\ZZ}{\ensuremath{\mathbb{Z}}}

\usepackage{amssymb,a4wide,epsfig,psfrag}
\usepackage{url}
\newtheorem{theorem}{Theorem}
\newtheorem{proposition}[theorem]{Proposition}
\newtheorem{lemma}[theorem]{Lemma}

\newfont{\Ar}{msam10}

\DeclareMathOperator{\Syl}{Syl}
\newcommand{\RR}{\ensuremath{\mathbb{R}}}
\newcommand{\SnWR}{S^n_{\rm wr}}

\DeclareMathOperator{\Aut}{Aut}

\DeclareMathOperator{\SL}{SL}
\DeclareMathOperator{\Sp}{Sp}
\DeclareMathOperator{\GL}{GL}
\DeclareMathOperator{\PSL}{PSL}
\DeclareMathOperator{\PGL}{PGL}

\DeclareMathOperator{\Min}{Min}

\title{On the integral homology of $\PSL_4(\ZZ)$ and other arithmetic groups}

\author{Mathieu Dutour Sikiri\'c}
\address{M. Dutour Sikiri\'c, Rudjer Boskovi\'c Institute, Bijenicka 54, 10000 Zagreb, Croatia}
\email{mdsikir@irb.hr}

\author{Graham Ellis}
\address{G. Ellis, Mathematics Department, National University of Ireland, Galway}
\email{graham.ellis@nuigalway.ie}

\author{Achill Sch\"urmann}
\address{A. Sch\"urmann, Institute of Mathematics, 
University of Rostock, 
18051 Rostock, 
Germany}
\email{achill.schuermann@uni-rostock.de}

\thanks{First author has been supported by Marie Curie fellowship MTKD-CT-2006-042685 and by the Croatian Ministry of Science, Education and Sport under contract 098-0982705-2707.}

\begin{document}

\begin{abstract}
We determine the integral homology of $\PSL_4(\ZZ)$ in degrees $\le 5$ and determine its $p$-part in higher degrees for the primes $p\ge 5$.
Our method applies to other arithmetic groups; as illustrations we include descriptions of the integral homology of $\PGL_3(\ZZ[i])$ and $\PGL_3(\ZZ[\exp(2\pi{i}/3)])$ in degrees $\le 5$.
\end{abstract}

\maketitle

\section{Introduction}
We obtain the following partial description of the integral homology of the projective special linear group $\PSL_4(\ZZ)$
(writing $H_{(p)}$ to denote the $p$-primary part of an Abelian group $H$).

\begin{theorem}\label{ONE}~

\begin{enumerate}
\item[(i)] $H_n(\PSL_4(\ZZ),\ZZ)=
\left\{ \begin{array}{ll}
 0 & n=1 \\
(\ZZ_2)^3 &n=2\\
\ZZ \oplus (\ZZ_4)^2 \oplus (\ZZ_3)^2 \oplus \ZZ_5, &n=3\\
(\ZZ_2)^4 \oplus \ZZ_5, &n=4\\
(\ZZ_2)^{13}, &n=5.
\end{array} \right. $\medskip
\item[(ii)] $H_n(\PSL_4(\ZZ),\ZZ)_{(5)}=
\left\{ \begin{array}{ll}
(\ZZ_5)^2 & n\equiv 3 {\rm\ mod\ } 4 \ (n\ge 6) \\
0 &n\equiv 0, 1, 2 {\rm\ mod\ } 4 \ (n\ge 6)
\end{array} \right. $\medskip
 \item[(iii)] $H_n(\PSL_4(\ZZ),\ZZ)_{(p)}=0$ for all primes $p\ge 7$ and $n\ge 0$.\medskip
\item[(iv)] $H_n(\PSL_4(\ZZ),\ZZ)$ is finite for all positive $n\ne 3$.
\end{enumerate}
\end{theorem}

The proof of Theorem \ref{ONE}(i) involves four steps.
The first, described in Section \ref{sectwo}, uses computer calculations
of perfect forms to explicitly determine a
 CW-structure on a homotopy retract $X$ of the space $S^4_{> 0}$ of $4\times 4$ positive definite symmetric matrices.
The retract $X$ is $6$-dimensional, contractible and admits a cellular action
of $G=\PSL_4(\ZZ)$ in which each cell $e$ has a finite stabilizer group $G^e$.
(We remark that an analogous $3$-dimensional retract of $S^3_{> 0}$, due to
Avner Ash, was used by Christophe Soul\'e \cite{soule1978} in his calculation
of the integral cohomology of $\SL_3(\ZZ)$.
Furthermore, versions of our retract $X$ for $\PSL_4(\ZZ)$ have been determined by Ash et al. in \cite{gunnells} and by Elbaz-Vincent et al. in \cite{quelques,parallel}, who used it to compute cohomology rationally and integrally at large primes only.) 
We will refer to \cite{schurmannBook} and \cite{martinet} for full details
on perfect forms  and  limit our exposition to those facts necessary for
Theorem \ref{ONE}.

The next step is described in Section \ref{secthree} and uses algorithms
from \cite{jsc,ellis} to determine small free $\ZZ G^e$-resolutions
of $\ZZ$ for each stabilizer group. Implementations of
these algorithms are available in the HAP \cite{hap} software package.
However, Theorem \ref{ONE} was obtained using an implementation
specially adapted by the first author.

The third step is described in Section \ref{secfour}.
It uses a generalization of a lemma of C.T.C. Wall to combine the
resolutions for the stabilizer groups with the cellular chain
complex $C_\ast(X)$ to form a free $\ZZ G$-resolution $R_\ast$ of $\ZZ$.
(For a full exposition of this technique we will refer to \cite{elliswilliams}
where it was illustrated with the detailed hand calculation of free
resolutions for generalized triangle groups acting on hyperbolic $3$-space.)
The technique has been automated in the software package \cite{hap}
and such an automation was used to obtain six terms of a free
$\ZZ G$-resolution $R_{\ast}$.

 The final step in the proof of Theorem \ref{ONE}(i) is routine: the Smith Normal Form algorithm is used to compute the homology of the chain complex $R_\ast\otimes_{\ZZ G}\ZZ$.

\medskip
Theorem \ref{ONE}(ii)-(iv) are standard applications of the Leray spectral sequence. Details are given in Section \ref{secfive}. The following result, which
 gives bounds on the $3$-primary part of the homology, is also derived in Section \ref{secfive}. 

\begin{proposition}\label{TWO} At the prime $p=3$ the Leray spectral sequence
$$E_{p,q}^1 =\bigoplus_{[e^p] } H_q(G^{e^p},\ZZ^\rho)_{(3)} \Rightarrow H_n(\PSL_4(\ZZ),\ZZ)_{(3)}$$
(where $[e^p]$ ranges over the $G$-orbits of $p$-dimensional $e^p$ cells in X and $\rho$ is a suitable action of $G^{e^p}$ on $\ZZ$) has first page $$E^1_{p,q}=
\left\{ \begin{array}{llll}
(\ZZ_3)^{q/2} &~~p=0,q=2+4k\ & \\
(\ZZ_3)^{(5+q)/2} &~~p=0,q=3+4k\  & \\
(\ZZ_3)^{(q+2)/4} &~~p=1,q=2+4k\ &  \\
(\ZZ_3)^{(q+5)/4} &~~p=1,q=3+4k\  & \\
\ZZ_3             &~~p=3,q=1+2k\         &\\
\ZZ_3             &~~p=4,q=1+4k\   \\
(\ZZ_3)^{(q+2)/4} &~~p=4,q=2+4k\   \\
(\ZZ_3)^{(q+5)/4} &~~p=4,q=3+4k\   \\
%0                       &~~p=4,q=4+4k\   \\
\ZZ_3             &~~p=5,q=1+2k\   \\
\ZZ_3             &~~p=6,q=1+4k\   \\
0                       &~~{\rm otherwise}
\end{array} (k\ge 0).\right. $$
\end{proposition}

The techniques underlying Theorem~\ref{ONE} can in principle be applied to other arithmetic groups. For example, 
they yield the following descriptions of the low-dimensional integral
 homology of the 3-dimensional projective general linear groups 
$\PGL_3(\ZZ[i])$, $\PGL_3(\ZZ[\omega])$ over the Gaussian  and Eisenstein integers (where $\omega = \exp(2\pi {i}/3)$).

\begin{theorem}\label{THREE}~\nopagebreak

$H_n(\PGL_3(\ZZ[i]),\ZZ)=
\left\{ \begin{array}{ll}
 0 & n=1 \\
(\ZZ_2)^2 &n=2\\
\ZZ \oplus \ZZ_4 \oplus \ZZ_8 \oplus \ZZ_3, &n=3\\
\ZZ \oplus \ZZ_2, &n=4\\
(\ZZ_2)^{5} \oplus (\ZZ_4)^{2}, &n=5.
\end{array} \right.$
\medskip

$H_n(\PGL_3(\ZZ[\omega]),\ZZ)=
\left\{ \begin{array}{ll}
 \ZZ_3 & n=1 \\
\ZZ_3 &n=2\\
\ZZ \oplus  \ZZ_8 \oplus (\ZZ_3)^3, &n=3\\
\ZZ \oplus (\ZZ_3)^2, &n=4\\
(\ZZ_2)^{2} \oplus (\ZZ_3)^{3}, &n=5.
\end{array} \right.$
\end{theorem}

The method of Theorem~\ref{ONE}(i) and Theorem~\ref{THREE} has also
been used to obtain the first few terms of an explicit resolution for
the symplectic group $\Sp_4(\ZZ)$. Although the integral (co)homology
has previously been calculated in this case by other means 
(see~\cite{brownsteinANDlee}), 
an advantage of an explicit resolution is that in principle
it can be used to compute the cohomology and explicit cocycles for finite index
subgroups such as congruence subgroups. See the final section for details.

We will not give details for the proof of Theorem \ref{THREE}.
Instead, in Section \ref{secsix} we explain how the results of the theorem
can be computed using the HAP \cite{hap} software package.
(We remark that the integral homology of $\PSL_2(\mathbb {O})$ for various
rings of quadratic integers $\mathbb {O}$ have been calculated by
J. Schwermer and K. Vogtmann in \cite{vogtmann} and by A. Rahm
and M. Fuchs \cite{rahm,rahm2}.)

\section{Perfect forms and a  cell complex with
  $\PSL_4(\ZZ)$-action.}\label{sectwo}

In the space $S^n_{>0}$ of positive definite, real symmetric $n\times n$ matrices,
we consider a specific $\SL_n(\ZZ)$ invariant polyhedral complex. 
Such complexes are classically studied in the arithmetic theory of
quadratic forms. One of them goes back to the  
work of Voronoi \cite{voronoiPerfect} on {\em perfect quadratic forms}.

For $A\in S^n_{>0}$ an associated positive definite quadratic form on $\RR^n$
is defined by $A[x]=x^t A x$. By this correspondence we simply
identify quadratic forms and symmetric matrices.
The {\em arithmetical minimum} of $A\in S^n_{>0}$ is defined by
$$
\min(A)=\min_{0\not= v\in \ZZ^n} A[v]
.
$$
The finite set of its {\em representatives} is denoted by
\begin{equation*}
\Min(A)=\{v \in \ZZ^n\mbox{~with~} A[v] = \min(A)\}.
\end{equation*}
We consider the set of forms having a fixed arithmetical
minimum, say~$1$:
$$
S^{n}_{=1}
=
\{ A \in S^n_{>0} \mbox{~with~} \min(A) = 1\}.
$$

This set is the boundary of a locally finite polyhedral set
known as {\em Ryshkov polyhedron} (see~\cite{schurmannBook} for details).
In particular, $S^{n}_{=1}$
is a piecewise linear surface of
dimension~${n + 1\choose 2} -1$ and the support of
a cell complex obtained from the faces of the Ryshkov polyhedron.  Moreover, $S^{n}_{=1}$ is contractible.
Each $v\in \ZZ^n$
determines one of the top-dimensional cells by the linear condition
$A[v]=1$. The $0$-dimensional cells in $S^{n}_{=1}$ correspond to 
perfect quadratic forms. These are characterized as being uniquely determined
by their arithmetical minimum (here~$1$) and its representatives.

The group $\SL_n(\ZZ)$ acts on $S^{n}_{>0}$ and $S^{n}_{=1}$ by
$A \mapsto P A P^t$. Some higher dimensional cells of $S^{n}_{=1}$
(and in particular all top-dimensional ones) 
have infinite stabilizer groups.
To avoid problems arising from such stabilizers,
we consider a deformation retract. 
A quadratic form $A\in S^n_{>0}$ is said to be {\em well-rounded} if there exist
linearly independent $v_1, \dots, v_n \in \Min (A)$.
The set $\SnWR$  of well rounded forms in $S^n_{=1}$ defines an
$n \choose 2$-dimensional polyhedral subcomplex, in which all cells 
are bounded polyhedra and thus have finite stabilizer groups. 
By a result of Ash \cite{ash}, $\SnWR$ is a homotopy retract of $S^n_{>0}$ (and $S^n_{=1}$)
and is thus contractible. In fact, $\SnWR$ is a minimal deformation
retract as 
there is no proper closed subset of $\SnWR$ which is  $\SL_n(\ZZ)$
invariant and contractible (see~\cite{pettet}).

The $0$-dimensional cells of $\SnWR$ (the perfect forms)
can be enumerated up to $\SL_n(\ZZ)$ equivalence using 
{\em Voronoi's algorithm}. In essence it is a graph traversal 
search on the $1$-dimensional subcomplex of $S^{n}_{=1}$ 
(see \cite{voronoiPerfect,martinet,schurmannBook} for more details).
Having a complete list of inequivalent perfect forms one can obtain
all other orbits of cells, as these  all contain at least one of
the perfect forms.

In dimension $4$ there are just two perfect forms up to $\SL_4(\ZZ)$
equivalence (associated to the root lattices
$\mathsf{A}_4$ and $\mathsf{D}_4$). 
These yield a $\PSL_4(\ZZ)$-equivariant CW-decomposition of $\SnWR$  involving the following cells and cell stabilizer groups:

\begin{itemize}
\item[] Dim. $0$:  two cell orbits with stabilizers $A_5$ and $(A_4 \times A_4) : C_2$.
\item[] Dim. $1$:  two cell orbits with stabilizers $S_3$ and $S_3 \times S_3$.
\item[] Dim. $2$:  two cell orbits, both with stabilizer $C_2\times C_2$.
\item[] Dim. $3$:  four cell orbits, two with stabilizer $C_2\times C_2$ and two with stabilizer $S_4$.
\item[] Dim. $4$:  four cell orbits with  stabilizers $S_3$,
$D_{8}$, $S_4$ and $C_2\times S_3\times S_3$.
\item[] Dim. $5$:  three cell orbits with  stabilizers $D_{24}$, $S_4$ and $A_5$.
\item[] Dim. $6$:  one cell orbit with stabilizer $((C_2\times C_2\times C_2\times C_2) : C_3) : C_2$.
\end{itemize}

\section{Resolutions for stabilizer subgroups}\label{secthree}
The method outlined in Section \ref{secfour} for computing $n$ terms of a free $\ZZ G$-resolution for $G=\PSL_4(\ZZ)$ 
will require $n-p$ terms of a free $\ZZ G^e$-resolution for each stabilizer group
of a $p$-dimensional cell in $\SnWR$. 
 An algorithm for computing a reasonably small 
free $\ZZ G^e$-resolution $R^{G^e}_\ast \rightarrow \ZZ$ for a
 generic finite group $G^e$ was described in \cite{jsc}. The   
implementation of this algorithm in \cite{jsc} produces, for instance,  a
free $\ZZ A_5$-resolution with $26$ free generators in dimension 6. The algorithm could be used to obtain $6-p$ terms for each of the stabilizer groups in the above tessellation of $\SnWR$.   

However, for computational efficiency it is necessary to craft
smaller resolutions for some  of the
stabilizer groups $G^e$. 
One way to do this is  to find a suitable  polytope $P$ on which $G^e$ acts 
 with small cell stabilizers. Then,  using Lemma \ref{wall}  
 below, we can  combine resolutions for those  subgroups of $G^e$ stabilizing cells of $P$ 
 with the cellular chain complex $C_\ast(P)$ in order to produce a free $\ZZ G^e$-resolution $R^{G^e}_\ast \rightarrow \ZZ$. The idea was illustrated in 
\cite{elliswilliams} with the hand  calculation of a free $\ZZ A_4$-resolution 
 involving $n+1$ free generators in degree $n$. The same technique can be used, for instance, to construct a free $\ZZ A_5$-resolution involving $n+1$ free generators in degree $n$; in this case we take $P$ to be the icosahedron and let $A_5$ act on $P$ as a subgroup of the symmetry group $C_2\times A_5$ of the icosahedron.    

The first author used this polytopal method to craft small resolutions for the larger stabilizer groups in the tessellation of
$\SnWR$.

\medskip
Both, the algorithm in \cite{jsc} and the polytopal method
\cite{ellis,elliswilliams} can return a free $\ZZ G^e$-resolution
$R^{G^e}_\ast$ endowed with a contracting homotopy.
Latter is encoded as a sequence of $\ZZ$-linear homomorphisms $h_n\colon R_n^{G^e} \rightarrow R_{n+1}^{G^e}$ satisfying
\begin{equation*}
h_{n-1}d_{n}+d_{n+1}h_n=1\ \ (n\ge0, h_{-1}=0).
\end{equation*}
This contracting homotopy is required by Lemma \ref{wall} below.

\section{Perturbation theory}\label{secfour}
It was observed in \cite{ellis} that a method of Wall \cite{wall} for 
constructing  free resolutions for  group extensions can be extended to a method for  constructing 
 free resolutions for groups acting  on contractible cellular spaces with nice stabilizer subgroups. We recall the method.  

Suppose a chain complex 
\[C_\ast \colon \cdots\rightarrow C_p\rightarrow C_{p-1} \rightarrow \cdots \rightarrow C_0 \rightarrow \ZZ ,\]
of $\ZZ G$-modules is given, where the modules $C_p$ are not necessarily $\ZZ G$-free.  Suppose also  
 that for each $p\ge 0$ we are given  a free $\ZZ G$-resolution $A_{p,\ast}$ of
the module $C_p$ 
\[A_{p,\ast} \colon \rightarrow A_{p,q} \rightarrow A_{p,q-1} \rightarrow
\cdots \rightarrow A_{p,0} \rightarrow  C_p \]
endowed with a contracting homotopy 
$h_q\colon A_{p,q}\rightarrow A_{p,q+1}$ ($q\ge 0$).

The following   lemma 
 explains how to construct a 
  chain complex $R_\ast$ of free $\ZZ G$-modules $R_n =\bigoplus_{p+q=n} A_{p,q}$ and a surjective  $\ZZ G$-chain map 
$$\phi_\ast\colon R_\ast \rightarrow C_\ast$$
 which induces  homology isomorphisms $H_n(\phi_\ast)\colon H_n(R_\ast) \stackrel{\cong}{\rightarrow} H_n(C_\ast)$ for $n\ge 0$.
 We state a slightly more general version of the lemma than was given in \cite{ellis}. However, the proof is the same and boils down to the proof of the
 special case given in \cite{wall}.

\begin{lemma}\label{wall}\cite{wall,ellis}
\begin{enumerate}
\item[(i)] 
Let $A_{p,q}$ ($p,q\ge 0$) be a bigraded family of free 
$\ZZ G$-modules. Suppose that there are $\ZZ G$-module homomorphisms $d_0\colon A_{p,q} \rightarrow
A_{p,q-1}$ such that $(A_{p,\ast},d_0)$ is an acyclic chain complex for each $p$. Set $C_p=H_0(A_{p,\ast},d_0)$ and suppose further that there are
$\ZZ G$-homomorphisms $\delta\colon C_p\rightarrow C_{p-1}$ for which
$(C_\ast,\delta)$ is a chain complex. Then there exist $\ZZ G$-homomorphisms $d_k\colon A_{p,q} \rightarrow A_{p-1,q+k-1}$ ($k\ge 1, p>k$) such that
$$d=d_0 + d_1 + \cdots \colon R_n=\bigoplus_{p+q=n}A_{p,q} \rightarrow R_{n-1}=\bigoplus_{p+q=n-1} A_{p,q}$$
is a differential on a chain complex $R_\ast$ of free $\ZZ G$-modules.
\item[(ii)] The canonical chain maps $\phi_p\colon A_{p,\ast} \rightarrow H_0(A_{p,\ast})$ constitute a chain map $\phi_\ast\colon R_\ast \rightarrow C_\ast$ which is an homology isomorphism. 
\item[(iii)]
 Suppose that there exist $\ZZ$-homomorphisms $h_0\colon A_{p,q}\rightarrow A_{p,q+1}$ such that $d_0h_0d_0(x) =d_0(x)$
for all $x\in A_{p,q+1}$. Then we can construct
 $d_k$ by first lifting $\delta$ to $d_1\colon A_{p,0}\rightarrow A_{p-1,0}$ and recursively defining $d_k=-h_0(\sum_{i=1}^k d_id_{k-i})$ on free generators of the module $A_{p,q}$.
\item[(iv)]
 Suppose that $C_\ast$ is acyclic, that  $H_0(S_\ast)\cong \ZZ$ and that
each $C_p$ is a free $\ZZ$-module. We can construct $\ZZ$-module homomorphisms
$h\colon R_n \rightarrow R_{n+1}$ satisfying $dhd(x)=d(x)$ by setting
$h(a_{p,q})=h_0(a_{p,q})-hd^+h_0(a_{p,q})+\epsilon(a_{p,q})$
for free generators $a_{p,q}$ of the module $A_{p,q}$.
Here $d^+=\sum_{i=1}^pd_i$ and, for $q\ge 1$, $\epsilon=0$. For $q=0$ we define $\epsilon=h_1-h_0d^+h_1+hd^+h_1+hd^+h_0d^+h_1$
where $h_1\colon A_{p,0}\rightarrow A_{p+1,0}$ is a $\ZZ$-linear homomorphism induced by a contracting homotopy on $C_\ast$.
\end{enumerate}
\end{lemma}

For $G=\PSL_4(\ZZ)$ and $X=\SnWR$ the cellular chain complex $C_\ast(X)$ can be viewed as a complex of $\ZZ G$-modules. Moreover, 
$$C_p(X) \cong \bigoplus_{[e^p]}\ \ \ZZ G \otimes_{\ZZ G^{e^p}} \ZZ $$
where $[e^p]$ ranges over the orbits of $p$-dimensional cells. A free $\ZZ G$-resolution $A_{p,\ast}$ of $C_p(X)$ can thus be obtained as
$$A_{p,\ast} \cong \bigoplus_{[e^p]}\ \  R_\ast^{G^{e^p}} \otimes_{\ZZ G^{e^p}} \ZZ $$
where  $R_\ast^{G^{e^p}}$ is any free $\ZZ G^{e^p}$-resolution of $\ZZ$. Moreover, contracting homotopies on $R_\ast^{G^{e^p}}$ induce a contracting homotopy on $A_{p,\ast}$.
 Lemma \ref{wall} thus gives an automated procedure to combine the chain complex $C_\ast(X)$ with resolutions for the stabilizer groups to produce a free $\ZZ G$-resolution $R_\ast^G$ of $\ZZ$.

\medskip
The first six terms of such a resolution $R_\ast^G$ were computed and, via the Smith Normal Form algorithm, used to prove Theorem \ref{ONE}(i).

\section{The Leray spectral sequence} \label{secfive}
Details on the  Leray spectral sequence can be found in \cite{brown}. 
For $G=\PSL_4(\ZZ)$ and $X=\SnWR$ it has the form 
$$E_{p,q}^1 =\bigoplus_{[e^p] } H_q(G^{e^p},\ZZ^\rho) \Rightarrow H_n(\PSL_4(\ZZ),\ZZ)$$
where $[e^p]$ ranges over the $G$-orbits of $p$-dimensional cells $e^p$ in the tessellation of $X$ and the integer coefficients are twisted by an action $\rho\colon G^{e^p}\rightarrow \Aut(\ZZ)=\{\pm 1\}$.

\medskip
Each stabilizer group is finite.  The homology of any finite group is finite in degrees $\ge 1$. So since $X$ is $6$-dimensional
 the spectral sequence immediately implies that $H_n(G,\ZZ)$ is finite for $n\ge 6$. A computer calculation shows that the
 homomorphism $E_{6,0}=\ZZ \rightarrow E_{5,0}=\ZZ^3$ is injective, and  $H_6(G,\ZZ)$ is thus finite.  Theorem \ref{ONE}(iv) therefore follows from Theorem~\ref{ONE}(i).

\medskip
No prime $p\ge 7$ divides the order of
 any of the stabilizer subgroups $G^{e}$ and thus no $p\ge 7$ divides the order of $H_q(G^e,\ZZ^\rho)$. So Theorem \ref{ONE}(iii) follows from the spectral sequence.

\medskip
The only $3$-groups and $5$-groups arising as subgroups of stabilizer groups 
are $C_3$,  $C_3\times C_3$ and $C_5$. These prime-power groups can only act trivially on $\ZZ$. Using Cartan and Eilenberg's 
identification of the $p$-part $H^q(G^e,\ZZ^\rho)_{(p)}$ of the cohomology  of $G^e$ with 
the $G$-stable elements in the cohomology 
$H^q(\Syl_p(G^e),\ZZ)$ of the Sylow $p$-subgroup, it is 
straightforward  to determine 
$H_{q-1}(G^e,\ZZ^\rho)_{(p)} = H^q(G^e,\ZZ^\rho)_{(p)}$ for $p=3,5$.
The $3$-part is presented in Proposition \ref{TWO}.
The $5$-part is given in the following.

\begin{proposition}\label{FIVE} At the prime $p=5$ the  Leray spectral sequence
$E_{p,q}^1  \Rightarrow H_n(\PSL_4(\ZZ),\ZZ)_{(5)}$
 has first page $$E^1_{p,q}=
\left\{ \begin{array}{llll}
\ZZ_5             &~~p=0,q=3+4k\     & \\
\ZZ_5             &~~p=5,q=3+4k\   \\
0                       &~~{\rm otherwise}
\end{array} (k\ge 0).\right. $$
\end{proposition}

Proposition \ref{FIVE} directly implies that
 $$H_n(\PSL_4(\ZZ),\ZZ)_{(5)}=
\left\{ \begin{array}{ll}
 0 {\rm \ or\ } (Z_5)^2 & ~~{\rm if~} n\equiv  3 {\rm\ mod\ } 4 \ (n\ge 6) \\
0 & ~~{\rm if~} n\equiv 0, 1, 2   {\rm\ mod\ } 4 \ (n\ge 6) \end{array} \right. $$
The alternating group $A_5$ is the stabilizer group of a $0$-cell of $X$.    A computer calculation shows that the inclusion $A_5\hookrightarrow 
\PSL_4(\ZZ)$ induces a surjection in cohomology $H^4(\PSL_4(\ZZ),\ZZ)_{(5)} \rightarrow H^4(A_5,\ZZ)_{(5)}$. The ring $H^{\ast}(A_5,\ZZ)_{(5)}\cong \ZZ_5[x]$ is generated by a single class $x$ in degree $4$. Hence the
 ring homomorphism  $H^\ast(\PSL_4(\ZZ),\ZZ)_{(5)} \rightarrow H^\ast(A_5,\ZZ)_{(5)}$   is surjective and consequently
 $H_n(\PSL_4(\ZZ),\ZZ)_{(5)}\cong H^{n+1}(\PSL_4(\ZZ),\ZZ)_{(5)}$ is non-trivial 
for $n\equiv 3 {\rm \ mod\ } 4$ ($n \ge 7$). This proves Theorem \ref{ONE}(ii).

\section{The groups $\PGL_3(\ZZ[i])$ and $\PGL_3(\ZZ[\omega])$} \label{secsix}

The above techniques apply to other arithmetic groups. 
In particular, the first author using the construction
of \cite{schurmannPerfect} for equivariant perfect forms
has computed two $8$-dimensional contractible
CW-spaces on which the groups $\PGL_3(\ZZ[i])$ and $\PGL_3(\ZZ[\omega])$
act cellularly with finite stabilizer groups.
These spaces, together with the above $6$-dimensional space
for $\PSL_4(\ZZ)$, have been stored in the HAP  package \cite{hap}
for the GAP computer algebra system.

The following short GAP session illustrates how the space for
$\PGL_3(\ZZ[i])$ can be accessed and used to compute:
(i) five terms of a free resolution and
(ii) the fourth integral homology group.

\begin{verbatim}
gap> C:=ContractibleGcomplex("PGL(3,Z[i])");;
gap> R:=FreeGResolution(C,5);;
gap> Homology(TensorWithIntegers(R),4);
[ 2, 0 ]
\end{verbatim}

The resolution can also be used to determine the homology of finite index subgroups of $\PGL_3(\ZZ[i])$, though GAP's standard implementation of the Smith Normal Form algorithm does not work well when the index is large.

%% PREVIOUS THEOREM, replaced by a remark in the introduction:
%\begin{theorem}\label{FOURTH}~\nopagebreak
%
%$H_n(\Sp_4(\ZZ),\ZZ)=
%\left\{\begin{array}{ll}
%\ZZ_2 & n=1 \\
%\ZZ_2 \oplus \ZZ &n=2\\
%(\ZZ_2)^2 \oplus \ZZ_4 \oplus \ZZ_{16} (\ZZ_3)^2\oplus \ZZ_5, &n=3\\
%\ZZ_2 \oplus \ZZ_4 \oplus \ZZ_3 , &n=4\\
%(\ZZ_2)^3 \oplus \ZZ_4 \oplus \ZZ_8 \oplus (\ZZ_3)^2 \oplus \ZZ_5, &n=5.
%\end{array} \right.$
%
%\end{theorem}

For the group $\Sp_4(\ZZ)$  the $4$-dimensional CW-complex
of \cite{Connell1993}, which comes also from the perfect forms in
dimension $4$, has been stored in the HAP package \cite{hap}. So,
 for  instance, the following GAP session computes
 the  homology $H_3(\Sp_4(\ZZ),\ZZ) = (\ZZ_2)^2\oplus \ZZ_{12} \oplus \ZZ_{240}$ in agreement with calculations of \cite{brownsteinANDlee}.

\begin{verbatim}
gap> C:=ContractibleGcomplex("Sp(4,Z)");;
gap> R:=FreeGResolution(C,4);;
gap> Homology(TensorWithIntegers(R),3);
[ 2, 2, 12, 240 ]
\end{verbatim}

\section{Acknowledgment}
We thank Mark McConnell and Paul Gunnells for help in using the cellular
decomposition of \cite{Connell1993}.

\end{document}